# Pseudo-maximization and self-normalized processes*

**Victor H. de la Peña[1], Michael J. Klass[2] and Tze Leung Lai[3]**

[1]*Department of Statistics, Columbia University*
*New York, NY 10027*
*e-mail:* vp@stat.columbia.edu

[2]*Department of Statisics, University of California*
*Berkeley, CA 94720-3860*
*e-mail:* klass@stat.berkeley.edu

[3]*Department of Statistics, Stanford University*
*Stanford, CA 94305-4065*
*e-mail:* lait@stat.stanford.edu

**Abstract:** Self-normalized processes are basic to many probabilistic and statistical studies. They arise naturally in the the study of stochastic integrals, martingale inequalities and limit theorems, likelihood-based methods in hypothesis testing and parameter estimation, and Studentized pivots and bootstrap-*t* methods for confidence intervals. In contrast to standard normalization, large values of the observations play a lesser role as they appear both in the numerator and its self-normalized denominator, thereby making the process scale invariant and contributing to its robustness. Herein we survey a number of results for self-normalized processes in the case of dependent variables and describe a key method called "pseudo-maximization" that has been used to derive these results. In the multivariate case, self-normalization consists of multiplying by the inverse of a positive definite matrix (instead of dividing by a positive random variable as in the scalar case) and is ubiquitous in statistical applications, examples of which are given.

**In Memory of Joseph L. Doob**

## 1. Introduction

This paper presents an introduction to the theory and applications of self-normalized processes in dependent variables, which was relatively unexplored until recently due to difficulties caused by the highly non-linear nature of self-normalization. We overcome these difficulties by using the method of mixtures which provides a tool for "pseudo-maximization".

*Research supported by NSF.

We dedicate this paper to the memory of J. L. Doob, a great probabilist who generously pointed one of us in this general direction, though we each had our independent initial path. While the first author was visiting the University of Illinois at Urbana–Champaign in the early 1990's, he took frequent hiking trips with Doob. Upon reaching a mountain-top in one of these trips, he asked Doob what were the most important open problems in probability. Doob replied that there were results in harmonic analysis involving harmonic functions divided by subharmonic or superharmonic functions that did not yet have analogues in the probabilistic setting of martingales. Guided by Doob's answer, de la Peña (1999) developed a new technique for obtaining exponential bounds for martingales. Subsequently, de la Peña, Klass and Lai (2000, 2004) introduced another method, which we call *pseudo-maximization*, to derive exponential and $L_p$-bounds for self-normalized processes. Via separate and newly innovated methods a universal LIL was also obtained. In this survey we review these results as well as those by others that fall (broadly) under the rubric of self-normalization. Our choice of topics has been guided by the usefulness and definitiveness of the results, and the light they shed on various aspects of probabilistic/statistical theory.

Self-normalized processes arise naturally in statistical applications. In standard form (as when connected to CLT's) they are unit-free and often permit the weakening or even the elimination of moment assumptions. The prototypical example of a self-normalized random process is Student's $t$-statistic. It replaces the population standard deviation $\sigma$ in $\sqrt{n}(\bar{X} - \mu)/\sigma$ by the sample standard deviation. Let $\{X_i\}$ be i.i.d normal $N(\mu, \sigma^2)$,

$$\bar{X}_n = \frac{\sum_{i=1}^n X_i}{n} \quad s_n^2 = \frac{\sum_{i=1}^n (X_i - \bar{X}_n)^2}{n-1}.$$

Then $T_n = \frac{\bar{X}_n - \mu}{s_n/\sqrt{n}}$ has the $t_{n-1}$-distribution. Let $Y_i = X_i - \mu$, $A_n = \sum_{i=1}^n Y_i$, $B_n^2 = \sum_{i=1}^n Y_i^2$. We can re-express $T_n$ in terms of $A_n/B_n$ as

$$T_n = \frac{A_n/B_n}{\sqrt{(n - (A_n/B_n)^2)/(n-1)}}.$$

Thus, information on the properties of $T_n$ can be derived from the self-normalized process

$$\frac{A_n}{B_n} = \frac{\sum_{i=1}^n Y_i}{\sqrt{\sum_{i=1}^n Y_i^2}}. \tag{1.1}$$

Hence, in a more general context, a self-normalized process assumes the form $A_n/B_n$, or $A_t/B_t$ in continuous time, where $B_t$ is a random variable that is used to estimate a dispersion measure of the process $A_t$.

Although the methodology of self-normalization dates back to 1908 when William Gosset (aka Student) introduced Student's $t$-statistic, active development of the probability theory of self-normalized processes began in the 1990's after the seminal work of Griffin and Kuelbs (1989, 1991) on laws of the iterated logarithm (LIL) for self-normalized sums of i.i.d. variables belonging to the domain of attraction of a normal or stable law. In particular, Bentkus and Götze (1996) derived a Berry-Esseen bound for Student's $t$-statistic, and Giné, Götze and Mason (1997) proved that the $t$-statistic has a limiting standard normal distribution if and only if $X_1$ is in the domain of attraction of a normal law, by making use of exponential and $L_p$ bounds for $A_n/B_n$, where $A_n = \sum_{i=1}^n X_i$ and $B_n^2 = \sum_{i=1}^n X_i^2$. The limiting distribution of the $t$-statistic when $X_1$ belongs to the domain of attraction of a stable law had been studied earlier by Efron (1969) and Logan, Mallows, Rice and Shepp (1973). Whereas Giné, Götze and Mason's result settles one of the conjectures of Logan, Mallows, Rice and Shepp that the self-normalized sum "is asymptotically normal if (and perhaps only if) $X_1$ is in the domain of attraction of the normal law," Chistyakov and Götze (2004) have settled that other conjecture that the "only possible nontrivial limiting distributions" are those when $X_1$ follows a stable law. Shao (1997) proved large deviation results for $\Sigma_{i=1}^n X_i / \sqrt{\Sigma_{i=1}^n X_i^2}$ without moment conditions and moderate deviation results when $X_1$ is the domain of attraction of a normal or stable law. Subsequently Shao (1997) obtained Cramér-type large deviation results when $E|X_1|^3 < \infty$. Jing, Shao and Zhou (2003) derived saddlepoint approximations for Student's $t$-statistic with no moment assumptions. Bercu, Gassiat and Rio (2002) obtained large and moderate deviation results for self-normalized empirical processes. Self-normalized sums of independent but non-identically distributed $X_i$ have been considered by Bentkus, Bloznelis and Götze (1996), Wang and Jing (1999) and Jing, Shao and Wang (2002). Chen (1999), Worms (2000) and Bercu (2001) have provided extensions to self-normalized sums of functions of ergodic Markov chains and autoregressive models. Giné and Mason (1998) relate the LIL of self-normalized sums of i.i.d. random variables $X_i$ to the stochastic boundedness of $(\Sigma_{i=1}^n X_i^2)^{-1/2} \Sigma_{i=1}^n X_i$.

Egorov (1998) gives exponential inequalities for a centered variant of (1.1). Along a related line of work, Caballero, Fernández and Nualart (1998) provide moment inequalities for a continuous martingale divided by its quadratic variation, and use these results to show that if $\{M_t, t \geq 0\}$ is a continuous martingale, null at zero, then for every $1 \leq p < q$, there exists a universal constant $C = C(p,q)$ such that

$$\|\frac{M_t}{\langle M \rangle_t}\|_p \leq C \|\frac{1}{\langle M \rangle_t^{\frac{1}{2}}}\|_q. \tag{1.2}$$

Related work in Revuz and Yor (1999 page 167) for continuous local martingales yields for all $p > q > 0$ the existence of a constant $C_{pq}$ such that

$$E \frac{(\sup_{s<\infty} |M_s|)^p}{\langle M \rangle_\infty^{q/2}} \leq C_{pq} E (\sup_{s<\infty} |M_s|)^{p-q}. \tag{1.3}$$

In Section 2 we describe the approaches of de la Peña (1999) to exponential inequalities for strong-law norming and of de la Peña, Klass and Lai (2000, 2004) to exponential inequalities and $L_p$-bounds of self-normalized processes.

Section 3 considers laws of the iterated logarithm for self-normalized martingales. Section 4 concludes with a discussion and review of self-normalization and pseudo-maximization in statistical applications.

## 2. Self-normalization and pseudo-maximization

### 2.1. Motivation

<div align="center">"BEHIND EVERY LIMIT THEOREM<br>THERE IS AN INEQUALITY"</div>

This folklore has been attributed to Kolmogorov.

**Example 2.1.** Let $\frac{S_n}{a_n}$ be a sequence of random variables. Then, to show that $\frac{S_n}{a_n} \to \mu$ in probability, Markov's inequality is often used:

$$P(|\frac{S_n}{a_n} - \mu| > \epsilon) \leq E|\frac{S_n}{a_n} - \mu|^p/\epsilon^p.$$

The weak law of large numbers for sums of i.i.d. random variables with finite variance uses the case $p = 2$ and the fact that the variance of the sum is the sum of the the variances. What happens when the variance is infinite, and when $a_n$ depends on $(X_1, X_2, \ldots, X_n)$ ?

**Example 2.2** (Almost Sure Growth Rate of a Sum). For non-decreasing $a_n$, if we can show for some $K > 0$ and for all $\epsilon > 0$ that

$$P\{\frac{S_n}{a_n} > (1 + 3\epsilon)K \text{ i.o.}\} = 0,$$

then we can conclude that $\limsup \frac{S_n}{a_n} \leq K$ a.s. By the Borell-Cantelli lemma, it suffices to show that for some $1 < n_1 < n_2 < \ldots$ with $a_j < (1 + \epsilon)a_{n_k}$ for $n_k \leq j < n_{k+1}$ whenever $k$ is sufficiently large,

$$\sum_{k=1}^{\infty} P(\max_{1 \leq j < n_{k+1}} S_j > (1 + \epsilon)Ka_{n_k}) < \infty.$$

Problems of this type frequently reduce to the use of Markov's inequality and finding appropriate bounds on $E \exp(\lambda S_n/a_n)$ to show that the above series converges. We are particularly interested in situations in which $a_n$ depends on the available data and hence is random. This further motivates the study of self-normalized sums.

**Example 2.3.** Consider the autoregressive model

$$Y_n = \alpha Y_{n-1} + \epsilon_n, \ Y_0 = 0,$$

where $\alpha$ is a fixed (unknown) parameter and $\epsilon_n$ are independent standard normal random variables. The maximum likelihood estimate (MLE) $\hat{\alpha}_n$ of $\alpha$ is the maximizer of the log likelihood function

$$\log f_\alpha(Y_1, ..., Y_n) = \sum_{j=1}^{n}(Y_j - \alpha Y_{j-1})^2/2 - n\log(\sqrt{2\pi}).$$

Differentiating with respect to $\alpha$ and equating to zero, we obtain

$$\hat{\alpha}_n = \frac{\sum_{j=1}^{n} Y_{j-1} Y_j}{\sum_{j=1}^{n} Y_{j-1}^2} = \frac{\sum_{j=1}^{n} Y_{j-1}(\alpha Y_{j-1} + \epsilon_j)}{\sum_{j=1}^{n} Y_{j-1}^2} = \alpha + \frac{\sum_{j=1}^{n} Y_{j-1}\epsilon_j}{\sum_{j=1}^{n} Y_{j-1}^2}.$$

Therefore, $\hat{\alpha}_n - \alpha$ can be expressed as a self-normalized random variable:

$$\hat{\alpha}_n - \alpha = \frac{\sum_{j=1}^{n} Y_{j-1}\epsilon_j}{\sum_{j=1}^{n} Y_{j-1}^2}. \tag{2.1}$$

In (2.1), the numerator $A_n := \sum_{j=1}^{n} Y_{j-1}\epsilon_j$ is a martingale with respect to the filtration $\mathcal{F}_n := \sigma(Y_1, ..., Y_n;\ \epsilon_1, ..., \epsilon_n)$. The denominator of (2.1) is

$$B_n^2 := \sum_{j=1}^{n} E[(Y_{j-1}\epsilon_j)^2|\mathcal{F}_{n-1}] = \sum_{j=1}^{n} Y_{j-1}^2,$$

which is the conditional variance of $A_n$. Therefore, $\hat{\alpha}_n - \alpha = A_n/B_n^2$ is a process self-normalized by the conditional variance. Since the $\epsilon_j$ are $N(0,1)$, for all $n \geq 1$ and $-\infty < \lambda < \infty$,

$$M_n := \exp\{\lambda \sum_{j=1}^{n} Y_{j-1}\epsilon_j - \frac{\lambda^2 \sum_{j=1}^{n} Y_{j-1}^2}{2}\}$$

is an exponential martingale and therefore satisfies the canonical assumption below, which we use to develop a comprehensive theory for self-normalized processes.

## 2.2. Canonical assumption and exponential bounds for strong law

We consider a pair of random variables $A$, $B$, with $B > 0$ such that

$$E \exp\{\lambda A - \lambda^2 B^2/2\} \leq 1. \tag{2.2}$$

There are three regimes of interest: when (2.2) holds

(i) for all real $\lambda$,
(ii) for all $\lambda \geq 0$,
(iii) for all $0 \leq \lambda < \lambda_0$, where $0 < \lambda_0 < \infty$.

Results are presented in all three cases. Such canonical assumptions imply various moment and exponential bounds, including the following bound connected to the law of large numbers in de la Peña (1999) for case (i).

**Theorem 2.1.** *Under the canonical assumption for all real $\lambda$,*

$$P(A/B^2 > x, 1/B^2 \leq y) \leq \exp(-\frac{x^2}{2y}) \qquad (2.3)$$

*for all $x, y > 0$.*

*Proof.* The key here is to "keep" the indicator when using Markov's inequality. In fact, for all measurable sets **S**,

$$P\{A/B^2 > x, \mathbf{S}\} = P\{\exp(A) > \exp(xB^2), \mathbf{S}\}$$
$$\leq \inf_{\lambda > 0} E[\exp\{\frac{\lambda}{2}A - \frac{\lambda}{2}xB^2\}1(A/B^2 > x, \mathbf{S})]$$
$$= \inf_{\lambda > 0} E[\exp\{\frac{\lambda}{2}A - \frac{\lambda^2}{4}B^2 - (\frac{\lambda}{2}x - \frac{\lambda^2}{4})B^2\}1(A/B^2 > x, \mathbf{S})]$$
$$\leq \sqrt{E\exp\{\lambda A - \frac{\lambda^2}{2}B^2\})}\sqrt{E[\exp\{-(\lambda x - \frac{\lambda^2}{2})B^2\}1(A/B^2 > x, \mathbf{S}))]},$$

by the Cauchy-Schwarz inequality. The first term in the last inequality is bounded by 1, by the canonical assumption. The value minimizing the second term is $\lambda = x$, and therefore

$$P(A/B^2 > x, S) \leq \sqrt{E[\exp\{\frac{-x^2 B^2}{2}\}1(A/B^2 > x, \mathbf{S})]}.$$

Letting $\mathbf{S} = \{1/B^2 \leq y\}$ gives the desired result. □

Applying this bound with $y = 1/z$ to both $A_n$ and $-A_n$ in Example 2.3 yields

$$P(|\hat{\alpha}_n - \alpha| > x, \sum_{j=1}^{n} Y_{j-1}^2 \geq z) \leq 2\exp(-x^2 z/2)$$

for all $x, z > 0$. The following variant of Theorem 2.1 generalizes a result of Lipster and Spokoiny (1999).

**Theorem 2.2.** *Under the canonical assumption (2.2) for all real $\lambda$,*

$$P(|A|/B > x, b \leq B \leq bs) \leq 4\sqrt{e}x(1 + \log s)\exp\{-x^2/2\},$$

*for all $b > 0$, $s \geq 1$ and $x \geq 1$.*

**Example 2.4** (Martingales and Supermartingales). The Appendix provides several classes of random variables that satisfy the canonical assumption (2.2) in a series of lemmas (A.1–A.8). These lemmas are closely related to martingale theory. Moreover, Lemmas A.2–A.8 are about a supermartingale condition (2.7) that is stronger than (2.2) for the regime $0 \leq \lambda < \lambda_0$.

Theorem 2.1 gives an inequality related to the law of large numbers. There is a class of results of this type in Khoshnevisan (1996), who points out that it essentially dates back to McKean (1962). A related reference is Freedman (1975).

**Theorem 2.3.** *Let $\{M_t, \mathcal{F}_t; 0 \leq t \leq \infty\}$ be a continuous martingale such that $E \exp(\gamma M_\infty) < \infty$ for all $\gamma$. Assume that $M_0 = 0$, that $\mathcal{F}_0$ contains all null sets and that the filtration $\{\mathcal{F}_t, 0 \leq t \leq \infty\}$ is right-continuous. Then, for any $\alpha, \beta, \lambda > 0$,*

$$P\{M_\infty \geq (\alpha + \beta \langle M \rangle_\infty)\lambda\} \leq \exp(-2\alpha\beta\lambda^2). \tag{2.4}$$

Khoshnevisan (1996) has also shown that if one assumes that the local modulus of continuity of $M_t$ is in some sense deterministic, then the inequality (2.4) can be reversed (up to a multiplicative constant). As applications of this result, he presents some large deviation results for such martingales. A related paper of Dembo (1996) provides moderate deviations for martingales with bounded jumps. Concerning moderate deviations, a general context for extending results like Theorem 2.3 to the case of discrete-time martingales can be found in de la Peña (1999), who provides a decoupling method for obtaining sharp extensions of exponential inequalities for martingales to the quotient of a martingale divided by its quadratic variation. In what follows we present three related results from de la Peña (1999). The first result is for martingales in continuous time and the last two involves discrete-time processes. The third result is a sharp extension of Bernstein's and Bennett's inequalities. (See also de la Peña and Giné (1999) for details.)

**Theorem 2.4.** *Let $\{M_t, \mathcal{F}_t, t \geq 0\}$, $M_0 = 0$, be a continuous martingale for which $\exp\{\lambda M_t - \lambda^2 \langle M \rangle_t / 2\}, t \geq 0$, is a supermartingale for all $\lambda > 0$. Let $\mathbf{A}$ be a set measurable with respect to $\mathcal{F}_\infty$. Then, for all $0 < t < \infty, \beta > 0, \alpha, x \geq 0$,*

$$P(M_t \geq (\alpha + \beta \langle M \rangle_t)x, \mathbf{A})$$
$$\leq E[\exp\{-x^2(\frac{\beta^2}{2}\langle M \rangle_t + \alpha\beta)\} \mid M_t \geq (\alpha + \beta \langle M \rangle_t)x, \mathbf{A}],$$

*and hence*

$$P\{M_t \geq (\alpha + \beta \langle M \rangle_t)x, \frac{1}{\langle M \rangle_t} \leq y \text{ for some } t < \infty\} \leq \exp\{-x^2(\frac{\beta^2}{2y} + \alpha\beta)\}.$$

**Theorem 2.5.** *Let $\{M_n = \sum_{i=1}^n d_i, \mathcal{F}_n, n \geq 0\}$ be a sum of conditionally symmetric random variables ($\mathcal{L}(d_i|\mathcal{F}_{i-1}) = \mathcal{L}(-d_i|\mathcal{F}_{i-1})$ for all $i$). Let $\mathbf{A}$ be a set* *measurable with respect to $\mathcal{F}_\infty$. Then, for all $n \geq 1$, $\beta > 0$, $\alpha$, $x \geq 0$,*

$$P(M_n \geq (\alpha + \beta \sum_{i=1}^n d_i^2)x, \mathbf{A})$$
$$\leq E[\exp\{-x^2(\frac{\beta^2}{2}\sum_{i=1}^n d_i^2 + \alpha\beta)\} \mid M_n \geq (\alpha + \beta \sum_{i=1}^n d_i^2)x, \mathbf{A}],$$

*and hence*

$$P\{M_n \geq (\alpha+\beta\sum_{i=1}^n d_i^2)x, \frac{1}{\sum_{i=1}^n d_i^2} \leq y \text{ for some } n \geq 1\} \leq \exp\{-x^2(\frac{\beta^2}{2y}+\alpha\beta)\}.$$

**Theorem 2.6.** *Let $\{M_n = \sum_{i=1}^n d_i, \mathcal{F}_n, n \geq 1\}$ be a martingale with $E(d_j^2|\mathcal{F}_{j-1}) = \sigma_j^2 < \infty$ and let $V_n^2 = \sum_{j=1}^n \sigma_j^2$. Furthermore assume that $E(|d_j|^k \mid \mathcal{F}_{j-1}) \leq (k!/2)\sigma_j^2 c^{k-2}$ a. s. for all $k > 2$ and some $c > 0$. Then, for all $\mathcal{F}_\infty$-measurable sets $\mathbf{A}$, $x > 0$,*

$$P\Big(\frac{M_n}{V_n^2} \geq x, \mathbf{A}\Big) \leq E\Big[\exp\Big\{-\Big(\frac{x^2}{1+cx+\sqrt{1+2cx}}\Big)V_n^2\Big\} \mid \frac{M_n}{V_n^2} \geq x, \mathbf{A}\Big],$$

*and hence*

$$P\Big(\frac{M_n}{V_n^2} \geq x, \frac{1}{V_n^2} \leq y \text{ for some } n\Big)$$
$$\leq \exp\Big\{-\frac{1}{y}\Big(\frac{x^2}{1+cx+\sqrt{1+2cx}}\Big)\Big\} \leq \exp\Big\{-\frac{x^2}{2y(1+cx)}\Big\}.$$

### 2.3. Pseudo-maximization (Method of Mixtures)

Note that if the integrand $\exp\{\lambda A - \lambda^2 B^2/2\}$ in (2.2) can be maximized over $\lambda$ inside the expectation (as can be done if $A/B^2$ is non-random), taking $\lambda = A/B^2$ would yield $E \exp(\frac{A^2}{2B^2}) \leq 1$. This in turn would give the optimal Chebyshev-type bound $P(\frac{A}{B} > x) \leq \exp(\frac{-x^2}{2})$. Since $A/B^2$ cannot (in general) be taken to be non-random, we need to find an alternative method for dealing with this maximization. One approach for attaining a similar effect involves integrating over a probability measure $F$, and using Fubini's theorem to interchange the order of integration with respect to $P$ and $F$. To be effective for all possible pairs $(A, B)$, the $F$ chosen would need to be as uniform as possible so as to include the maximum value of $\exp\{\lambda A - \lambda^2 B^2/2\}$ regardless of where it might occur. Thereby some mass is certain to be assigned to and near the random value $\lambda = A/B^2$ which maximizes $\lambda A - \lambda^2 B^2/2$. Since all uniform measures are multiples of Lebesgue measure (which is infinite), we construct a finite measure (or a sequence of finite measures) which tapers off to zero at $\lambda = \infty$ as slowly as we can manage. This approach will be used in what follows to provide exponential and moment inequalities for $A/B$, $A/\sqrt{B^2 + (EB)^2}$, $A/\{B\sqrt{\log\log(B \vee e^2)}\}$. We begin with the second case where the proof is more transparent. The approach, pioneered by Robbins and Siegmund (1970) and commonly known as the method of mixtures, was used by de la Peña, Klass and Lai (2004) to prove the following.

**Theorem 2.7.** *Let $A, B$ with $B > 0$ be random variables satisfying the canonical assumption (2.2) for all $\lambda \in \mathbf{R}$. Then*

$$P(\frac{|A|}{\sqrt{B^2 + (EB)^2}} > x) \leq \sqrt{2}\exp(-x^2/4) \qquad (2.5)$$

*for all $x > 0$.*

*Proof.* Multiplying both sides of (2.2) by $(2\pi)^{-1/2} y \exp(-\lambda^2 y^2/2)$ (with $y > 0$) and integrating over $\lambda$, we obtain by using Fubini's theorem that

$$\begin{aligned}
1 &\geq \int_{-\infty}^{\infty} E \frac{y}{\sqrt{2\pi}} \exp\left(\lambda A - \frac{\lambda^2}{2} B^2\right) \exp\left(-\frac{\lambda^2 y^2}{2}\right) d\lambda \\
&= E\Big[\frac{y}{\sqrt{B^2+y^2}} \exp\Big\{\frac{A^2}{2(B^2+y^2)}\Big\} \times \int_{-\infty}^{\infty} \frac{\sqrt{B^2+y^2}}{\sqrt{2\pi}} \\
&\quad \exp\{-\frac{B^2+y^2}{2}\Big(\lambda^2 - 2\frac{A}{B^2+y^2}\lambda + \frac{A^2}{(B^2+y^2)^2}\Big)\}d\lambda\Big] \\
&= E\Big[\frac{y}{\sqrt{B^2+y^2}} \exp\Big(\frac{A^2}{2(B^2+y^2)}\Big)\Big].
\end{aligned} \qquad (2.6)$$

By the Cauchy-Schwarz inequality and (2.6),

$$\begin{aligned}
&E \exp\Big\{\frac{A^2}{4(B^2+y^2)}\Big\} \\
&\leq \Big\{\Big(E\frac{y\exp\{\frac{A^2}{2(B^2+y^2)}\}}{\sqrt{B^2+y^2}}\Big)\Big(E\sqrt{\frac{B^2+y^2}{y^2}}\Big)\Big\}^{1/2} \\
&\leq \Big(E\sqrt{\frac{B^2}{y^2}+1}\Big)^{1/2}.
\end{aligned}$$

Since $E\sqrt{\frac{B^2}{y^2}+1} \leq E(\frac{B}{y}+1)$, the special case $y = EB$ gives

$$E\exp(A^2/[4(B^2 + (EB)^2)]) \leq \sqrt{2}.$$

Combining Markov's inequality with this yields

$$P\Big\{\frac{|A|}{\sqrt{B^2+(EB)^2}} \geq x\Big\} = P\Big\{\frac{A^2}{4(B^2+(EB)^2)} \geq \frac{x^2}{4}\Big\} \leq \sqrt{2}\,\exp(-x^2/4).$$

$\square$

In what follows we discuss the analysis of certain boundary crossing probabilities by using the method of mixtures, first introduced in Robbins and Siegmund (1970), under the following refinement of the canonical assumption:

$$\{\exp(\lambda A_t - \lambda^2 B_t^2/2),\ t \geq 0\} \text{ is a supermartingale with mean } \leq 1 \qquad (2.7)$$

for $0 \leq \lambda < \lambda_0$, with $A_0 = 0$. We begin by introducing the Robbins-Siegmund (R-S) boundaries; see Lai (1976). Let $F$ be a finite positive measure on $(0, \lambda_0)$ and assume that $F(0, \lambda_0) > 0$. Let $\Psi(u, v) = \int \exp(\lambda u - \lambda^2 v/2) dF(\lambda)$. Given $c > 0$ and $v > 0$, the equation

$$\Psi(u, v) = c$$

has a unique solution

$$u = \beta_F(v, c).$$

Moreover, $\beta_F(v, c)$ is a concave function of $v$ and

$$\lim_{v \to \infty} \frac{\beta_F(v, c)}{v} = b/2,$$

where

$$b = \sup\left\{y > 0 : \int_0^y F(d\lambda) = 0\right\},$$

with sup over the empty set equal to zero. The R-S boundaries $\beta_F(v, c)$ can be used to analyze the boundary crossing probability

$$P\{A_t \geq g(B_t) \text{ for some } t \geq 0\}$$

when $g(B_t) = \beta_F(B_t^2, c)$ for some $F$ and $c > 0$. This probability equals

$$\begin{aligned} &P\{A_t \geq \beta_F(B_t^2, c) \text{ for some } t \geq 0\} \\ =\ &P\{\Psi(A_t, B_t^2) \geq c \text{ for some } t \geq 0\} \leq F(0, \lambda_0)/c, \end{aligned} \quad (2.8)$$

applying Doob's inequality to the supermartingale $\Psi(A_t, B_t)$, $t \geq 0$.

**Example 2.5.** Let $\delta > 0$ and

$$dF(\lambda) = \frac{1}{\lambda(\log(1/\lambda))(\log\log(1/\lambda))^{1+\delta}} d\lambda$$

for $0 < \lambda < e^{-e}$. Let $\log_2(x) = \log\log(x)$ and $\log_3(x) = \log\log_2(x)$. As shown in Example 4 of Robbins and Siegmund (1970), for fixed $c > 0$,

$$\beta_F(v, c) = \sqrt{2v[\log_2 v + (3/2 + \delta)\log_3 v + \log(c/2\sqrt{\pi}) + o(1)]},$$

as $v \to \infty$. With this choice of $F$, the probability in (2.8) is bounded by $F(0, e^{-e})/c$ for all $c > 0$. Given $\epsilon > 0$, take $c$ large enough so that $F(0, e^{-e})/c < \epsilon$. Since $\epsilon$ can be arbitrarily small and since for fixed $c$, $\beta_F(v, c) \sim \sqrt{2v \log\log v}$ as $v \to \infty$,

$$\limsup \frac{A_t}{\sqrt{2B_t^2 \log\log B_t^2}} \leq 1,$$

on the set $\{\lim_{t \to \infty} B_t = \infty\}$.

### *2.4. Moment inequalities for self-normalized processes*

The inequality (1.2) obtained by Caballero, Fernández and Nualart (1998) is used by them to establish the continuity and uniqueness of the solutions of a nonlinear stochastic partial differential equation. A natural question that arises in connection with (1.2) is what about the case when the normalization is done by $\sqrt{\langle M \rangle_t}$. The following result from de la Peña, Klass and Lai (2000, 2004) provides an answer to this question.

**Theorem 2.8.** *Let $A, B > 0$ be two random variables satisfying (2.2) for all $\lambda > 0$. Let $\log^+(x) = 0 \vee \log x$. Then for all $p > 0$,*

$$E|\frac{A^+}{B}|^p \leq c_{1,p} + c_{2,p} E(\log^+ \log(B \vee \frac{1}{B}))^{p/2}, \qquad (2.9)$$

$$E(\frac{A^+}{B\sqrt{1 \vee \log^+ \log(B \vee \frac{1}{B})}})^p \leq C_p, \qquad (2.10)$$

*for some universal constants $0 < c_{1,p}, c_{2,p}, C_p < \infty$.*

The following example shows the optimality properties of Theorem 2.8.

**Example 2.6.** Let $\{Y_i\}$, $i = 1, \ldots$ be a sequence of i.i.d. Bernoulli random variables such that $P\{Y_i = 1\} = P\{Y_i = -1\} = \frac{1}{2}$. Let

$$T = \inf\{n \geq e^e : \sum_{i=1}^n Y_i \geq \sqrt{2n \log \log n}\},$$

with $T = \infty$ if no such $n$ exists. From the delicate LIL by Erdős (1942), it follows that $P(T < \infty) = 1$. Let $d_{n,j} = Y_j 1(T \geq j)$ if $1 \leq j \leq n$ and $d_{n,j} = 0$ if $j > n$. Then almost surely,

$$\frac{d_{n,1} + \ldots + d_{n,n}}{\sqrt{d_{n,1}^2 + \ldots + d_{n,n}^2}} = \frac{\sum_{j=1}^{\min\{T,n\}} Y_j}{\sqrt{\min\{T, n\}}} \to \sqrt{2 \log \log T},$$

as $n \to \infty$. Therefore, the second term in (2.9) can not be removed.

The proof of Theorem 2.8 is based on the following result. Let $L : (0, \infty) \to (0, \infty)$ be a nondecreasing function such that

$$L(cy) \leq 3cL(y) \text{ for all } c \geq 1 \text{ and } y > 0, \qquad (2.11a)$$

$$L(y^2) \leq 3L(y) \text{ for all } y \geq 1, \qquad (2.11b)$$

$$\int_1^\infty \frac{dx}{xL(x)} = \frac{1}{2}. \qquad (2.11c)$$

Then

$$f(\lambda) = \frac{1}{\lambda L(\max(\lambda, 1/\lambda))}, \quad \lambda > 0,$$

is a density. An example satisfying (2.11a,b,c) is

$$L(x) = 2(\log xe^e)(\log\log xe^e)^2 1(x \geq 1).$$

Let $g(x) = \frac{\exp\{x^2/2\}}{x} 1(x \geq 1)$. Then

$$E \frac{g(\frac{A}{B})}{L(\frac{A}{B}) \vee L(B \vee \frac{1}{B})} \leq \frac{3}{\int_0^1 \exp^{-x^2/2} dx}. \tag{2.12}$$

Making use of this, de la Peña, Klass and Lai (2000, 2004) extend the inequalities in (2.9) and (2.10) to $Eh\left(\frac{A^+}{B}\right)$ and $Eh\left(\frac{A^+}{B\sqrt{\log\log(e^e \vee B \vee \frac{1}{B})}}\right)$ for nonnegative non-decreasing functions $h(\cdot)$ such that $h(x) \leq g^\delta(x)$ for some $0 < \delta < 1$.

### 2.5. An expectation version of the LIL for self-normalized processes

We next study the case of self-normalized inequalities when there is no possibility of explosion at the origin by shifting the denominator away from zero. An important result in this direction comes from Graversen and Peskir (2000).

**Theorem 2.9.** *Let $\{M_t, \mathcal{F}_t, t \geq 0\}$ be a continuous local martingale with quadratic variation process $\langle M \rangle_t$, $t \geq 0$. Let $l(x) = \sqrt{\log(1 + \log(1 + x))}$. Then there exist universal constants $D_1$, $D_2 > 0$ such that*

$$D_1 E l(\langle M \rangle_\tau) \leq E\left(\frac{\sup_{0 \leq t \leq \tau} |M_s|}{\sqrt{1 + \langle M \rangle_t}}\right) \leq D_2 E l(\langle M \rangle_\tau),$$

*for all stopping times $\tau$ of $M$.*

The proof of this result was obtained by making use of Lenglart's inequality, the optional sampling theorem and Ito's formula. Shortly after this result appeared, de la Peña, Klass and Lai (2000) introduced the moment bounds in the last section, in which the denominator is not shifted away from 0. They then realized that shifted moment bounds can also be obtained for the case in which (2.2) or (2.7) only hold for $0 < \lambda < \lambda_0$. Subsequently, de la Peña, Klass and Lai (2004) proved part (i) of the following theorem for more general processes than continuous local martingales. Part (ii) of the theorem can be proved by arguments similar to those of Theorem 2 of de la Peña, Klass and Lai (2000).

**Theorem 2.10.** *Let $T = \{0, 1, 2, \ldots\}$ or $T = [0, \infty)$, $1 < q \leq 2$, and $\Phi_q(\theta) = \theta^q/q$. Let $A_t$ and $B_t > 0$ be stochastic processes satisfying*

$$\{\exp(\lambda A_t - \Phi_q(\lambda B_t)), t \in T\} \text{ is a supermartingale with mean} \leq 1 \tag{2.13}$$

*for $0 < \lambda < \lambda_0$ and such that $A_0 = 0$ and $B_t$ is nondecreasing in $t > 0$. In the case $T = [0, \infty)$, assume furthermore that $A_t$ and $B_t$ are right-continuous.*

Let $L : [1, \infty) \to (0, \infty)$ be a nondecreasing function satisfying (2.11a, b, c). Let $\eta > 0, \lambda_0 \eta > \delta > 0$, and $h : [0, \infty) \to [0, \infty)$ be a nondecreasing function such that $h(x) \leq e^{\delta x}$ for all large $x$.

(i) There exists a constant $\kappa$ depending only on $\lambda_0, \eta, \delta, q, h$ and $L$ such that

$$Eh\Big(\sup_{t \geq 0} \big\{A_t(B_t \vee \eta)^{-1}[1 \vee \log^+ L(B_t \vee \eta)]^{-(q-1)/q}\big\}\Big) \leq \kappa. \qquad (2.14)$$

(ii) There exists $\widetilde{\kappa}$ such that for any stopping time $\tau$,

$$Eh\Big(\sup_{0 \leq t \leq \tau} |(B_t \vee \eta)^{-1} A_t|\Big) \leq \widetilde{\kappa} + Eh(|\widetilde{\kappa} \vee \log^+ L(B_t \vee \eta)]^{(q-1)/q}). \qquad (2.15)$$

## 3. Self-normalized LIL for stochastic sequences

### 3.1. Stout's LIL: Self-normalizing via conditional variance

A well-known result in martingale theory is Stout's (1970, 1973) LIL that uses the square root of the conditional variance for normalization. The key to the proof of Stout's result is an exponential supermartingale in Lemma A.5 of the Appendix.

**Theorem 3.1.** *Let $\{d_n, \mathcal{F}_n\}, n = 1, \ldots$ be an adapted sequence with $E(d_n|\mathcal{F}_{n-1}) \leq 0$. Set $M_n = \sum_{i=1}^n d_i$, $\sigma_n^2 = \sum_{i=1}^n E(d_i^2|\mathcal{F}_{i-1})$. Assume that*
  (i) $d_n \leq m_n$ for $\mathcal{F}_{n-1}$-measurable $m_n \geq 0$,
  (ii) $\sigma_n^2 < \infty$ a.s for all $n$,
  (iii) $\lim_{n \to \infty} \sigma_n^2 = \infty$ a.s.,
  (iv) $\limsup_{n \to \infty} m_n \sqrt{\log \log(\sigma_n^2)}/\sigma_n = 0$ a.s.
*Then*
$$\limsup \frac{M_n}{\sqrt{2\sigma_n^2 \log \log \sigma_n}} \leq 1 \text{ a.s.}$$

### 3.2. Discrete-time martingales: self-normalizing via sums of squares

By making use of Lemma A.8 (see Appendix), de la Peña, Klass and Lai (2004) have proved the following upper LIL for self-normalized and suitably centered sums of random variables, under no assumptions on their joint distributions.

**Theorem 3.2.** *Let $X_n$ be random variables adapted to a filtration $\{\mathcal{F}_n\}$. Let $S_n = X_1 + \ldots + X_n$ and $V_n^2 = X_1^2 + \ldots + X_n^2$. Then, given any $\lambda > 0$, there exist positive constants $c_\lambda$ and $b_\lambda$ such that $\lim_{\lambda \to 0} b_\lambda = \sqrt{2}$ and*

$$\limsup \frac{S_n - \sum_{i=1}^n \mu_i[-\lambda v_n, c_\lambda v_n)}{V_n(\log \log V_n)^{1/2}} \leq b_\lambda \text{ a.s. on } \{\lim V_n = \infty\}, \qquad (3.1)$$

*where $v_n = V_n(\log \log V_n)^{-1/2}$ and $\mu_i[c, d) = E\{X_i 1(c \leq X < d)|\mathcal{F}_{i-1}\}$ for $c < d$.*

The constant $b_\lambda$ in Theorem 3.2 can be determined as follows. For $\lambda > 0$, let $h(\lambda)$ be the positive solution of $h - \log(1+h) = \lambda^2$. Let $b_\lambda = h(\lambda)/\lambda$, $\gamma = h(\lambda)/\{1 + h(\lambda)\}$, and $c_\lambda$ is determined via $\lambda/\gamma$. Then $\lim_{\lambda \to 0} b_\lambda = \sqrt{2}$. Let $e_k = \exp(k/\log k)$. The basic idea underlying the proof of Theorem 3.2 pertains to upper-bounding the probability of an event of the form $E_k = \{t_{k-1} \leq \tau_k < t_k\}$, in which $t_j$ and $\tau_j$ are stopping times defined by

$$
\begin{aligned}
t_j &= \inf\{n : V_n \geq e_j\}, \\
\tau_j &= \inf\{n \geq t_j : S_n - \textstyle\sum_{i=1}^n \mu_i[-\lambda v_n, c_\lambda v_n) \geq (1+3\epsilon)b_\lambda V_n(\log \log V_n)^{1/2}\},
\end{aligned}
\tag{3.2}
$$

letting $\inf \emptyset = \infty$. Note that for $i < n$, the centering constants $\mu_i[-\lambda v_n, c_\lambda v_n)$ involve $v_n$ that is not determined until time $n$, so the centered sums that result do not usually form a martingale. However, sandwiching $\tau_k$ between $t_{k-1}$ and $t_k$ enables us to replace both the random exceedance and truncation levels in (3.2) by constants. Then the event $E_k$ can be re-expressed in terms of two simultaneous inequalities, one involving centered sums and the other involving a sum of squares. These inequalities combine to permit application of Lemma A.8 (see Appendix). Thereby we conclude that $P(E_k)$ can be upper-bounded by the probability of an event involving the supremum of a nonnegative supermartingale with mean $\leq 1$, to which Doob's maximal inequality can be applied; see de la Peña, Klass and Lai (2004, pages 1924-1925) for details.

Although Theorem 3.2 gives an upper LIL for any adapted sequence of random variables $X_n$, the upper bound in (3.1) may not be attained. Example 6.4 of de la Peña, Klass and Lai (2004) suggests that one way to sharpen the bound is to center $X_n$ at its conditional median before applying Theorem 3.2 to $\widetilde{X}_n = X_n - \operatorname{med}(X_n | \mathcal{F}_{n-1})$. On the other hand, if $X_n$ is a martingale difference sequence such that $|X_n| \leq m_n$ a.s. for some $\mathcal{F}_{n-1}$-measurable $m_n$ with $m_n = o(v_n)$ and $v_n \to \infty$ a.s., then Theorem 6.1 of de la Peña, Klass and Lai (2004) shows that (3.1) is sharp in the sense that

$$
\limsup \frac{S_n}{V_n(\log \log V_n)^{1/2}} = \sqrt{2} \text{ a.s.} \tag{3.3}
$$

The following example of de la Peña, Klass and Lai (2004) illustrates the difference between self-normalization by $V_n$ and by $\sigma_n$ in Theorems 3.2 and 3.1.

**Example 3.1.** Let $X_1 = X_2 = 0$, $X_3, X_4, \ldots$ be independent random variables such that

$$P\{X_n = -1/\sqrt{n}\} = 1/2 - n^{-1/2}(\log n)^{1/2} - n^{-1}(\log n)^{-2},$$

$$P\{X_n = -m_n\} = n^{-1}(\log n)^{-2}, \ P\{X_n = 1/\sqrt{n}\} = 1/2 + n^{-1/2}(\log n)^{1/2}$$

for $n \geq 3$, where $m_n \sim 2(\log n)^{5/2}$ is chosen so that $EX_n = 0$. Then $P\{X_n = -m_n \text{ i.o.}\} = 0$. Hence, with probability 1, $V_n^2 = \Sigma_{i=1}^n i^{-1} + O(1) = \log n + O(1)$ and it is shown in de la Peña, Klass and Lai (2004, page 1926) that

$$
\frac{S_n}{V_n(\log \log V_n)^{1/2}} \sim \frac{4(\log n)^{3/2}}{3\{(\log n)(\log \log \log n)\}^{1/2}} \to \infty \text{ a.s.,} \tag{3.4}
$$

and that

$$\limsup \left( S_n - \sum_{i=1}^n EX_i 1(|X_i| \leq 1) \right) \Big/ \{V_n (\log \log V_n)^{1/2}\} = \sqrt{2} \text{ a.s.} \quad (3.5)$$

Note that $m_n/v_n \to \infty$ a.s. This shows that without requiring $m_n = o(v_n)$, the left-hand side of (3.3) may be infinite. On the other hand, $X_n$ is clearly bounded above and $\sigma_n^2 = \sum_{i=1}^n \text{Var}(X_i) \sim 4 \sum_{i=1}^n (\log i)^3/i \sim (\log n)^4$, yielding

$$\frac{S_n}{\sigma_n (\log \log s_n)^{1/2}} \sim \frac{4(\log n)^{3/2}}{3(\log n)^2 (\log \log \log n)^{1/2}} \to 0 \text{ a.s.,} \quad (3.6)$$

which is consistent with Stout's (1973) upper LIL. Note, however, that $m_n/\sigma_n \sim 2(\log n)^{1/2} \to \infty$ and therefore one still does not have Stout's (1970) lower LIL in this example.

## 4. Statistical applications

Most of the probability theory of self-normalized processes developed in the last two decades is concerned with $A_t$ self-normalized by $B_t$ in the case $A_t = \Sigma_{i \leq t} X_i$ is a sum of i.i.d. random vectors $X_i$ and $B_t = (\Sigma_{i \leq t} X_i X_i')^{1/2}$, using a key property that

$$E \exp\{\theta' X_1 - \rho(\theta' X_1)^2\} < \infty \text{ for all } \theta \in \mathbf{R}^k \text{ and } \rho > 0, \quad (4.1)$$

as observed by Chan and Lai (2000, pages 1646–1648). In the i.i.d. case, the finiteness of the moment generating function (4.1) enables one to embed the underlying distribution in an exponential family and one can then use change of measures (exponential tilting) to derive saddlepoint approximations for large or moderate deviation probabilities of self-normalized sums or for more general boundary crossing probabilities. Specifically, letting $C_n = \Sigma_{i=1}^n X_i X_i'$ and $\psi(\theta, \rho)$ denote the left hand side of (4.1), we have

$$E[\exp\{\theta' A_n - \rho \theta' C_n \theta - n\psi(\theta, \rho)\}] = 1. \quad (4.2)$$

Let $P_{\theta,\rho}$ be the probability measure under which $(X_i, X_i X_i')$ are i.i.d. with density function $\exp(\theta' X_i - \rho \theta' X_i X_i' \theta)$ with respect to $P = P_{0,0}$. The random variable inside the square brackets of (4.2) is simply the likelihood ratio statistic based on $(X_1, \ldots, X_n)$, or the Radon-Nikodym derivative $dP_{\theta,\rho}^{(n)}/dP^{(n)}$, where the superscript $(n)$ denotes restriction of the measure to the $\sigma$-field generated by $\{X_1, \ldots, X_n\}$. For the case of dependent random vectors, although we no longer have the simple cumulant generating function $n\psi(\theta, \rho)$ in (4.2) to develop precise large or moderate deviation approximations, we can still derive exponential bounds by applying the pseudo-maximization technique to (2.2) or (2.13), which is weaker than (4.2), as shown in the following two examples from de la Peña, Klass and Lai (2004, pages 1921-1922) who use a multivariate normal distribution with mean 0 and covariance matrix $V^{-1}$ for the mixing distribution $F$ to generalize (2.8) to the multivariate setting.

**Example 4.1.** Let $\{d_n\}$ be a sequence of random vectors adapted to a filtration $\{\mathcal{F}_n\}$ such that $d_i$ is conditionally symmetric. Then for any $a > 1$ and any positive definite $k \times k$ matrix $V$,

$$P\left\{\frac{(\Sigma_1^n d_i')(V + \Sigma_1^n d_i d_i')^{-1}(\Sigma_1^n d_i)}{\log\det(V + \Sigma_1^n d_i d_i') + 2\log(a/\sqrt{\det V})} \geq 1 \text{ for some } n \geq 1\right\} \leq \frac{1}{a}.$$

**Example 4.2.** Let $M_t$ be a continuous local martingale taking values in $\mathbf{R}^k$ such that $M_0 = 0$, $\lim_{t\to\infty} \lambda_{\min}(\langle M\rangle_t) = \infty$ a.s. and $E\exp(\theta'\langle M\rangle_t \theta) < \infty$ for all $\theta \in \mathbf{R}^k$ and $t > 0$. Then for any $a > 1$ and any positive definite $k \times k$ matrix $V$,

$$P\left\{\frac{M_t'(V + \langle M\rangle_t)^{-1} M_t}{\log\det(V + \langle M\rangle_t) + 2\log(a/\sqrt{\det\langle M\rangle_t})} \geq 1 \text{ for some } t \geq 0\right\} = \frac{1}{a}. \quad (4.3)$$

The reason why equality holds in (4.3) is that $\{\int f(\theta)\exp(\theta' M_t - \theta'\langle M\rangle_t \theta/2)d\theta, t \geq 0\}$ is a nonnegative continuous martingale to which an equality due to Robbins and Siegmund (1970) can be applied; see Corollary 4.3 of de la Peña, Klass and Lai (2004).

Self-normalization is ubiquitous in statistical applications, although $A_n$ and $C_n$ need no longer be linear functions of the observations $X_i$ and $X_i X_i'$ as in the $t$-statistic or Hotelling's $T^2$-statistic in the multivariate case. Section 4.1 gives an overview of self-normalization in statistical applications, and Section 4.2 discusses the connections of the pseudo-maximization approach with likelihood and Bayesian inference.

### 4.1. Self-normalization in statistical applications

The $t$-statistic $\sqrt{n}(\bar{X}_n - \mu)/s_n$ is a special case of more general *Studentized statistics* $(\widehat{\theta}_n - \theta)/\widehat{\mathrm{se}}_n$ that are of fundamental importance in statistical inference on an unknown parameter $\theta$ of an underlying distribution from which the sample observations $X_1, \ldots, X_n$ are drawn. In nonparametric inference, $\theta$ is a functional $g(F)$ of the underlying distribution function $F$ and $\widehat{\theta}_n$ is usually chosen to be $g(\widehat{F}_n)$, where $\widehat{F}_n$ is the empirical distribution. The standard deviation of $\widehat{\theta}_n$ is often called its *standard error*, which is typically unknown, and $\widehat{\mathrm{se}}_n$ denotes a consistent estimate of the standard error of $\widehat{\theta}_n$. For the $t$-statistic, $\mu$ is the mean of $F$ and $\bar{X}_n$ is the mean of $\widehat{F}_n$. Since $\mathrm{Var}(\bar{X}_n) = \mathrm{Var}(X_1)/n$, we estimate the standard error of $\bar{X}_n$ by $s_n/\sqrt{n}$, where $s_n^2$ is the sample variance. An important property of a Studentized statistic is that it is an *approximate pivot*, which means that its distribution is approximately the same for all $\theta$; see Efron and Tibshirani (1993, Section 12.5) who make use of this pivotal property to construct bootstrap-$t$ confidence intervals and tests. For parametric problems, $\theta$ is usually a multidimensional vector and $\widehat{\theta}_n$ is an asymptotically normal estimate (e.g., by maximum likelihood). Moreover, the asymptotic covariance matrix $\Sigma_n(\theta)$ of $\widehat{\theta}_n$ depends on the unknown parameter $\theta$, so $\Sigma_n^{-1/2}(\widehat{\theta}_n)(\widehat{\theta}_n - \theta)$ is the self-normalized (Studentized) statistic that can be used an approximate pivot for tests and confidence regions.

The theoretical basis for the approximate pivotal property of Studentized statistics lies in the limiting standard normal distribution, or in some other limiting distribution that does not involve $\theta$ (or $F$ in the nonparametric case). To derive the asymptotic normality of $\widehat{\theta}_n$, one often uses a martingale $M_n$ associated with the data, and approximates $\Sigma_n^{-1/2}(\widehat{\theta}_n)(\widehat{\theta}_n - \theta)$ by $\langle M \rangle_n^{-1/2} M_n$. For example, in the asymptotic theory of the maximum likelihood estimator $\widehat{\theta}_n$, $\Sigma_n(\theta)$ is the inverse of the observed Fisher information matrix $I_n(\theta)$, and the asymptotic normality of $\widehat{\theta}_n$ follows by using Taylor's theorem to derive

$$-I_n(\theta)(\widehat{\theta}_n - \theta) \doteq \sum_{i=1}^{n} \nabla \log f_\theta(X_i | X_1, \ldots, X_{i-1}), \tag{4.4}$$

The right-hand side of (4.4) is a martingale whose predictable variation is $-I_n(\theta)$. Therefore the Studentized statistic associated with the maximum likelihood estimator can be approximated by a self-normalized martingale.

### 4.2. Pseudo-maximization in likelihood and Bayesian inference

Let $X_1, \ldots, X_n$ be observations from a distribution with joint density function $f_\theta(x_1, \ldots, x_n)$. Likelihood inference is based on the likelihood function $L_n(\theta) = f_\theta(X_1, \ldots, X_n)$, whose maximization leads to the maximum likelihood estimator $\widehat{\theta}_n$. Bayesian inference is based on the posterior distribution of $\theta$, which is the conditional distribution of $\theta$ given $X_1, \ldots, X_n$ when $\theta$ is assumed to have a prior distribution with density function $\pi$. Under squared error loss, the Bayes estimator $\widetilde{\theta}_n$ is the mean of the posterior distribution whose density function is proportional to $L_n(\theta)\pi(\theta)$, i.e.,

$$\widetilde{\theta}_n = \int \theta \pi(\theta) L_n(\theta) d\theta \Big/ \int \pi(\theta) L_n(\theta) d\theta. \tag{4.5}$$

Recall that $L_n(\theta)$ is maximized at the maximum likelihood estimator $\widehat{\theta}_n$. Applying Laplace's asymptotic formula to the integrals in (4.5) shows that $\widetilde{\theta}_n$ is asymptotically equivalent to $\widehat{\theta}_n$, so integrating $\theta$ over the posterior distribution in (4.5) amounts to pseudo-maximization.

Let $\theta_0$ denote the true parameter value. A fundamental quantity in likelihood inference is the likelihood ratio martingale

$$\frac{L_n(\theta)}{L_n(\theta_0)} = e^{\ell_n(\theta)}, \text{ where } \ell_n(\theta) = \sum_{i=1}^{n} \log \frac{f_\theta(X_i | X_1, \ldots, X_{i-1})}{f_{\theta_0}(X_i | X_1, \ldots, X_{i-1})}. \tag{4.6}$$

Note that $\nabla \ell_n(\theta_0)$ is also a martingale; it is the martingale in the right-hand side of (4.4). Clearly $\int e^{\ell_n(\theta)} \pi(\theta) d\theta$ is also a martingale for any probability density function $\pi$ of $\theta$. Lai (2004) shows how the pseudo-maximization approach can be applied to $e^{\ell_n(\theta)}$ to derive boundary crossing probabilities for the generalized likelihood ratio statistics $\ell_n(\widehat{\theta}_n)$ that lead to efficient procedures in sequential analysis.

## Appendix

**Lemma A.1.** *Let $W_t$ be a standard Brownian motion. Assume that $T$ is a stopping time such that $T < \infty$ a.s.. Then*

$$E \exp\{\lambda W_T - \lambda^2 T/2\} \leq 1,$$

*for all $\lambda \in \mathbf{R}$.*

**Lemma A.2.** *Let $M_t$ be a continuous, square-integrable martingale, with $M_0 = 0$. Then $\exp\{\lambda M_t - \lambda^2 \langle M \rangle_t / 2\}$ is a supermartingale for all $\lambda \in \mathbf{R}$, and therefore*

$$E \exp\{\lambda M_t - \lambda^2 \langle M \rangle_t / 2\} \leq 1.$$

*If $M_t$ is only assumed to be a continuous local martingale, then the inequality is also valid (by application of Fatou's lemma).*

**Lemma A.3.** *Let $\{M_t : t \geq 0\}$ be a locally square-integrable martingale, with $M_0 = 0$. Let $\{V_t\}$ be an increasing process, which is adapted, purely discontinuous and locally integrable; let $V^{(p)}$ be its dual predictable projection. Set $X_t = M_t + V_t$, $C_t = \sum_{s \leq t}((\Delta X_s)^+)^2$, $D_t = \{\sum_{s \leq t}((\Delta X_s)^-)^2\}_t^{(p)}$, $H_t = \langle M \rangle_t^c + C_t + D_t$. Then $\exp\{X_t - V_t^{(p)} - \frac{1}{2}H_t\}$ is a supermartingale and hence*

$$E \exp\{\lambda(X_t - V_t^{(p)}) - \lambda^2 H_t / 2\} \leq 1 \text{ for all } \lambda \in \mathbf{R}.$$

Lemma A.2 follows from Proposition 3.5.12 of Karatzas and Shreve (1991). Lemma A.3 is taken from Proposition 4.2.1 of Barlow, Jacka and Yor (1986). The following lemma holds without any integrability conditions on the variables involved. It is a generalization of the fact that if $X$ is any symmetric random variable, then $A = X$ and $B = |X|$ satisfy the canonical condition (2.2); see de la Peña (1999).

**Lemma A.4.** *Let $\{d_i\}$ be a sequence of variables adapted to an increasing sequence of $\sigma$-fields $\{\mathcal{F}_i\}$. Assume that the $d_i$'s are conditionally symmetric (i.e., $\mathcal{L}(d_i|\mathcal{F}_{i-1}) = \mathcal{L}(-d_i|\mathcal{F}_{i-1})$). Then $\exp\{\lambda \Sigma_{i=1}^n d_i - \lambda^2 \Sigma_{i=1}^n d_i^2 / 2\}$, $n \geq 1$, is a supermartingale with mean $\leq 1$, for all $\lambda \in \mathbf{R}$.*

Note that any sequence of real-valued random variables $X_i$ can be "symmetrized" to produce an exponential supermartingale by introducing random variables $X_i'$ such that

$$\mathcal{L}(X_n'|X_1, X_1', \ldots, X_{n-1}, X_{n-1}', X_n) = \mathcal{L}(X_n|X_1, \ldots, X_{n-1})$$

and setting $d_n = X_n - X'_n$; see Section 6.1 of de la Peña and Giné (1999). For the next four lemmas, Lemma A.5 is the tool used by Stout (1970, 1973) to obtain the LIL for martingales, Lemma A.6 is de la Peña's (1999) extension of Bernstein's inequality for sums of independent variables to martingales, and Lemmas A.7 and A.8 correspond to Lemma 3.9(ii) and Corollary 5.3 of de la Peña, Klass and Lai (2004).

**Lemma A.5.** *Let $\{d_n\}$ be a sequence of random variables adapted to an increasing sequence of $\sigma$-fields $\{\mathcal{F}_n\}$ such that $E(d_n|\mathcal{F}_{n-1}) \leq 0$ and $d_n \leq M$ a.s. for all $n$ and some nonrandom positive constant $M$. Let $0 < \lambda_0 \leq M^{-1}$, $A_n = \sum_{i=1}^n d_i$, $B_n^2 = (1 + \frac{1}{2}\lambda_0 M)\sum_{i=1}^n E(d_i^2|\mathcal{F}_{i-1})$, $A_0 = B_0 = 0$. Then $\{\exp(\lambda A_n - \frac{1}{2}\lambda^2 B_n^2), \mathcal{F}_n, n \geq 0\}$ is a supermartingale for every $0 \leq \lambda \leq \lambda_0$.*

**Lemma A.6.** *Let $\{d_n\}$ be a sequence of random variables adapted to an increasing sequence of $\sigma$-fields $\{\mathcal{F}_n\}$ such that $E(d_n|\mathcal{F}_{n-1}) = 0$ and $\sigma_n^2 = E(d_n^2|\mathcal{F}_{n-1}) < \infty$. Assume that there exists a positive constant $M$ such that $E(|d_n|^k|\mathcal{F}_{n-1}) \leq (k!/2)\sigma_n^2 M^{k-2}$ a.s. or $P(|d_n| \leq M|\mathcal{F}_{n-1}) = 1$ a.s. for all $n \geq 1$, $k > 2$. Let $A_n = \sum_{i=1}^n d_i$, $V_n^2 = \sum_{i=1}^n E(d_i^2|\mathcal{F}_{i-1})$, $A_0 = V_0 = 0$. Then $\{\exp(\lambda A_n - \frac{1}{2(1-M\lambda)}\lambda^2 V_n^2), \mathcal{F}_n, n \geq 0\}$ is a supermartingale for every $0 \leq \lambda \leq 1/M$.*

**Lemma A.7.** *Let $\{d_n\}$ be a sequence of random variables adapted to an increasing sequence of $\sigma$-fields $\{\mathcal{F}_n\}$ such that $E(d_n|\mathcal{F}_{n-1}) \leq 0$ and $d_n \geq -M$ a.s. for all $n$ and some non-random positive constant $M$. Let $A_n = \sum_{i=1}^n d_i$, $B_n^2 = 2C_\gamma \sum_{i=1}^n d_i^2$, $A_0 = B_0 = 0$ where $C_\gamma = -\{\gamma + \log(1-\gamma)\}/\gamma^2$. Then $\{\exp(\lambda A_n - \frac{1}{2}\lambda^2 B_n^2, \mathcal{F}_n, n \geq 0\}$ is a supermartingale for every $0 \leq \lambda \leq \gamma M^{-1}$.*

**Lemma A.8.** *Let $\{\mathcal{F}_n\}$ be an increasing sequence of $\sigma$-fields and $y_n$ be $\mathcal{F}_n$-measurable random variables. Let $0 \leq \gamma_n < 1$ and $0 < \lambda_n \leq 1/C_{\gamma_n}$ be $\mathcal{F}_{n-1}$-measurable random variables, with $C_\gamma$ given in Lemma A.7. Let $\mu_n = E\{y_n 1(-\gamma_n \leq y_n < \lambda_n)|\mathcal{F}_{n-1}\}$. Then $\exp\{\sum_{i=1}^n (y_i - \mu_i - \lambda_i^{-1} y_i^2)\}$ is a supermartingale whose expectation is $\leq 1$.*